\documentclass{amsart}
\usepackage{cite}
\usepackage{amssymb}
\usepackage{color}  
\usepackage{mathtools}
\definecolor{mcpurple}{RGB}{70,0,200} 
\definecolor{mcblue}{RGB}{0,150,255} 
\definecolor{mcgreen}{RGB}{0,128,0}
\definecolor{mclime}{RGB}{0,200,0}
\definecolor{mcred}{RGB}{200,25,100}

\definecolor{mcbluegrey}{RGB}{225,230,240}
\definecolor{mcgreengrey}{RGB}{230,240,230}
\definecolor{mcmagenta}{RGB}{138,43,226}

\usepackage{mathrsfs}
\usepackage{amsmath}
\usepackage{thmtools, thm-restate}
\usepackage{verbatim}
\usepackage{hyperref} 
    \usepackage{cleveref}
\usepackage{memhfixc} 
\hypersetup{ 
    colorlinks=true,
    linkcolor=blue,    
    citecolor=blue,
    urlcolor=blue,
}
\usepackage{float}

\usepackage{tikz}
\usetikzlibrary{matrix}
\tikzset{diagram/.style={matrix of math nodes, inner sep=0pt, row
    sep=#1, column sep=2.5em, text height=1.5ex, text depth=.25ex,
    nodes={inner sep=1ex}}}
\tikzset{diagram/.default=2.5em}
\newcommand\diagram{\path node[diagram]}

\usepackage{memhfixc} 

\usepackage{amsthm}
\usepackage{thmtools, thm-restate}
\newtheorem{thm}{Theorem}[section]

\newtheorem{lem}[thm]{Lemma}
\newtheorem{cor}[thm]{Corollary}
\newtheorem{prop}[thm]{Proposition}
\newtheorem{defn}[thm]{Definition}

\newtheorem{rmk}[thm]{Remark}

\newcommand{\ZZ}{\mathbb{Z}} 
\newcommand{\QQ}{\mathbb{Q}}

\newcommand{\PP}{\mathscr{P}}

\newcommand{\OK}{\mathcal{O}_K}
\newcommand{\OL}{\mathcal{O}_L}

\newcommand{\pp}{\mathfrak{p}}

\newcommand{\qq}{\mathfrak{q}} 
\newcommand{\aaa}{\mathfrak{a}} 
\newcommand{\mmm}{\mathfrak{m}}

\newcommand{\Stab}{\operatorname{Stab}}
\newcommand{\surj}{\twoheadrightarrow}

\newcommand{\defeq }{\vcentcolon=}
\newcommand{\nequiv}{\not\equiv}
\newcommand{\ism}{\cong}

\newcommand{\inv}{^{-1}}

\newcommand{\vK}{\mathbf{x}_K}
\newcommand{\mKpos}{m_+}
\newcommand{\mKneg}{m_-}
\newcommand{\MK}{\mathbf{M}_4}
\newcommand{\MKG}{\mathbb{M}_{4,G}}

\newcommand{\spin}{\operatorname{spin}}
\newcommand{\Norm}{\operatorname{Norm}}
\newcommand{\ord}{\operatorname{ord}}
\newcommand{\Gal}{\operatorname{Gal}}

\newcommand{\Art}{\operatorname{Art}}

\newcommand{\preArt}[2]{\mathcal{A}_{#1}^{#2}}
\newcommand{\nRCF}[1]{\operatorname{nR}^{#1}}
\newcommand{\nRCG}[1]{{\operatorname{n\mathcal{C}l} }^{#1}}  %$\nRCG{\mm_0}$
\newcommand{\totpos}[1]{{#1}\succ 0}
\newcommand{\switchpf}[1]{#1}
\newcommand{\switch}[1]{}
\title {On the Asymptotics of a Prime Spin Relation, II}
\author {Christine McMeekin}
\date{November 2018}

\begin{document}

\maketitle

%%%%%%%%%%%%%%%%%%%%%%%%%%%%%%%%%%%
%\begin{center}
%\includegraphics[width=.4\textwidth]{APSRIIpicture.jpg}\\

%\small{ \textit{
%A cubic cyclic Galois group \\
%was blowing in the wind \\
%when we spotted a prime number\\
%and beside it\\
%a treasure map.\\
%We followed the treasure map\\
%to the golden units\\
%to find\\
%we were the stars\\
%all along.}}
%\end{center}
%%%%%%%%%%%%%%%%%%%%%%%%%%%%%%%

\begin{abstract}
For $K$ a cyclic cubic number field with odd class number containing a unit $\omega$ such that $\Norm_{K/\QQ}(\omega)=\Norm_{K/\QQ}(1-\omega)=-1$, we prove that the density of rational primes $p$ that satisfy the spin relation,
\[
\spin(\pp,\sigma) = \spin(\pp,\sigma\inv) \quad \text{for all } \sigma\neq 1 \in \Gal(K/\QQ),
\]
where $\pp$ is a prime of $K$ above $p$ is equal to $1/2$. 

Furthermore, we prove that this density restricted to rational primes that are $1 \bmod 4$ is $1/4$ and this density restricted to rational primes that are $-1 \bmod 4$ is $3/4$.
\end{abstract}

\tableofcontents

\section{Introduction}\label{sec:intro}%[thm:Dstarequality3]

%%%%%%%%%%%%%%%

For $K$ a cyclic number field of degree $\geq 3$ with odd class number $h(K)$ such that every totally positive unit is a square, the spin of an (odd) ideal is defined as follows.

\begin{defn}\cite{FIMR}\label{defn:spin}\switch{defn:spin} Let $\sigma\neq 1 \in \Gal(K/\QQ)$. Given an odd principal ideal $\aaa$, we define the \textit{spin} of $\aaa$ (with respect to $\sigma$) to be
\[
\spin(\aaa,\sigma) \defeq  \left(  \frac{\alpha}{ \aaa^\sigma}   \right),
\]
where $\aaa^{h(K)}=(\alpha)$, $\alpha$ is totally positive, and $\left( \frac{\cdot}{\cdot} \right)$ denotes the quadratic residue symbol in $K$. 
\end{defn}

%%%%%%%%%%%%%%%

\begin{restatable}{defn}{goldenunit}\label{defn:goldenunit}\switch{defn:goldenunit}
For $K$ a cyclic cubic number field with ring of integers $\OK$, if a unit $\omega\in \OK^\times$ satisfies
\[
\Norm(\omega) = \Norm(1 - \omega) = -1,
\]
then $\omega$ is said to be a \textit{golden unit}.
\end{restatable}

%%%%%%%%%%%%%%%

For $S\subseteq R$ sets of rational primes, let 
\[
d(R|S) = d\left( \frac{R}{S} \right)\defeq  \lim_{N\to \infty} \frac{\# R|_N}{\#S|_N}\]
where $S|_N$ and $R|_N$ denote the set of primes in $S$ and $R$ respectively of absolute norm less than $N\in \ZZ_+$.

%%%%%%%%%%%%%%%%%%%

%Let $K$ be a \switch{cubic cyclic} number field \switch{with odd class number}.
%%\todo{To Do: Check assumptions on $K$ needed to define $d_\pm$}
Let $S$ denote the set of odd rational primes that split completely in $K/\QQ$. 
Let 
\begin{align*}
&S_+\defeq \{ p \in S: p \equiv 1 \bmod 4\ZZ \}, \quad \text{and} \\
&S_-\defeq  \{ p \in S: p \equiv -1 \bmod 4\ZZ \}.
\end{align*}
For a fixed sign $\pm$, define
\[
 R_\pm \defeq \{p\in S_\pm: \spin(\pp,\sigma) = \spin(\pp,\sigma\inv)  \text{ for all } \sigma\neq 1 \in \Gal(K/\QQ) \},
 \]
 \[
 R\defeq  R_+ \cup R_-.
 \]

%\[ 
%d_\pm\defeq d(R_\pm|S_\pm).
%\]

%Define
%\[
%d_\pm\defeq  \lim_{N\to \infty} \frac{\# R_{\pm,N}}{\#S_{\pm,N}} \quad \text{and} \quad 
%d_K\defeq  \lim_{N\to \infty} \frac{\# R_{N}}{\#S_{N}}.
%\]

%%%%%%%%%%%%%%%

\begin{restatable}{thm}{Dstarequality}\label{thm:Dstarequality3}\switch{thm:Dstarequality3}
Let $K$ be a cyclic cubic number field with odd class number. 
Assume $K$ has a golden unit. %$\omega\in\OK^\times$
%and assume $5$ is inert in $K/\QQ$. 
Then $d(R|S)={1}/{2}$. More particularly,
\[
d(R_+|S_+) = \frac{1}{4}, \quad 
d(R_-|S_-) =\frac{3}{4}.
\]
\end{restatable}

%\begin{thm}\label{thm:Dstarequality3}\switch{thm:Dstarequality3}
%Let $K$ be a cyclic cubic number field with odd class number. Assume $K$ has a golden unit.Then $d_K=\frac{1}{2}$. More particularly,
%\[
%d_+ = \frac{1}{4}, \quad d_-=\frac{3}{4}.
%\]
%\end{thm}

%%%%%%%%%%%%%%

%%%%%%%%%%%%%%

%\todo{More general assumptions on K...}

Define $d_+:=d(R_+|S_+)$ and $d_-:=d(R_-|S_-)$.
The table below gives the computed positive and negative Starlight invariants and the corresponding restricted densities for the totally real cyclic number field of degree $n$ and conductor $\ell$.

Note that the density $d(R|S)$ is given by the average of $d_+$ and $d_-$ by Dirichlet/Chebotarev's density theorem where $R:=R_+\cup R_-$ and $S:=S_+\cup S_-$.

\begin{center}
\begin{tabular}{c c | c c | c c}
$n$ & $\ell$ & $\mKpos$ & $d_+$ & $\mKneg$ & $d_-$ \\ \hline
3  &  7  &  0  &  1/4  &  1  &  3/4 \\
5  &  11  &  0  &  1/16  &  1  &  5/16 \\
7  &  43  &  2  &  15/64  &  1  &  7/64 \\
11  &  23  &  0  &  1/1024  &  3  &  33/1024 \\
13  &  53  &  0  &  1/4096  &  5  &  65/4096 \\
\end{tabular}
\end{center}

The file ``Starlightpm.m" contains the code used to generate the table above.

%%%%%%%%%%%%%% treasuremap
\section{The Coincidence of Narrow and Wide}
%Class Field Theory When Narrow Equals Wide\\ and a Treasure Map %from Primes to Units
\label{sec:treasuremap:purple}%[treasuremap]%[RCFh]%[Uquadsub]

%%%%%%%%%%%%

%Let K tot real U_T=U^2...
%\begin{center}
 %   \switch{sec:treasuremap:purple}
%\end{center}

Let $K$ be a totally real number field,\switch{(cyclic? abelian?)} over $\QQ$ with odd class number $h(K)$. Assume 2 is inert in $K/\QQ$. Let $U\defeq \OK^\times$ denote the group of units of $K$, let $U_T$ denote the totally positive units, and let $U^2$ denote the square units.
When every totally positive unit of $K$ is a square unit (i.e. $U_T=U^2$), we can define a map from the prime ideals of $K$ to $U/U^2$. This map (defined in Lemma \ref{treasuremap}\switch{treasuremap}) will be used in Proposition  \ref{prop:prettypicture:existence}\switch{prop:prettypicture:existence}.
%to prove the existence of a non-trivial $\Gal(K/\QQ)$ orbit of $\MK$ satisfying $\star$.
 %which is one of two main ingredients to 
Proposition \ref{prop:prettypicture:existence}\switch{prop:prettypicture:existence} together with \ref{thm:pmdensityformulas}\switch{thm:pmdensityformulas} gives the main result,
Theorem \ref{thm:Dstarequality3}\switch{thm:Dstarequality3}.

%We now diverge momentarily from the spin of prime ideals to discuss some class field theory in the case when every totally positive unit is a square. 
We say a modulus $\mmm$ is \textit{narrow} whenever it is divisible by all infinite places. We say a modulus is \textit{wide} whenever it is not divisible by any infinite place. We say a ray class group or ray class field is narrow or wide whenever its defining modulus is narrow or wide respectively. 
For a modulus $\mmm$ of $K$, $\nRCF{\mmm}_K$ denotes the narrow ray class field over $K$ with conductor $\mmm\infty$ and  $\nRCG{\mmm}_K$ denotes the corresponding narrow ray class group, where $\infty$ denotes the product of all infinite places of $K$.

We restate Lemma 3.1 in \cite{APSR}. The proof is an exercise in class field theory.
    
    \begin{lem}\cite{APSR}\label{RCFh} \switch{RCFh}
Let $K$ be a totally real number field. The following are equivalent.
\begin{enumerate}
    \item $U_T=U^2$.
    \item The \textit{narrow} and \textit{wide} Hilbert class groups of $K$ coincide.
    \item Every principal ideal of $K$ has a totally positive generator.
\end{enumerate}

\end{lem}
%\switchpf{
%\begin{proof}
%To show the equivalence of (1) and (2), apply Theorem V.1.7 in \cite{MilneCFT} using the modulus given by the product of all infinite places. Statements (3) and (2) are equivalent by the definitions of narrow and wide class groups.
%\end{proof}
%}

%%%%%%%%%%%%

  \begin{lem}\label{Euu}\label{Uquadsub}\label{ExistenceKpp} \switch{Uquadsub} %%\todo{$K/\QQ$ cyclic? 2 inert? $n=[K:\QQ]$ odd? $n$ prime?}
   Let $K$ be a totally real number field, Galois over $\QQ$ such that the class number $h(K)$ is odd and ${U_T}={U^2}$. Let $\pp$ be an odd prime of $K$. Then the {narrow} ray class field over $K$ of conductor $\pp$ has a unique subextension that is quadratic over $K$. 
    \end{lem}
    \begin{proof}
    %%\todo{To Do: Check assumptions. $K/\QQ$ cyclic? 2 inert? $n=[K:\QQ]$ odd? $n$ prime? To Do: re-write this proof, defining notation from Milne.}
    
    We first show that the narrow ray class group over $K$ of conductor $\pp$ has even order. We then show the 2-part of the narrow ray class group over $K$ of conductor $\pp$ is cyclic.
    
    Let $\mmm$ be the narrow modulus with finite part $\pp$.
    Let 
    \begin{align*}
    &K_\mmm\defeq \{a\in K^\times: \ord_2(a)=0\},\\
    &K_{\mmm,1}\defeq \{a \in K^\times: \ord_2(a-1)\geq \ord_2(q), \totpos{a}\},\text{ and}\\
    &U_{\mmm,1}\defeq K_{\mmm,1}\cap U,
    \end{align*}
    where $\ord_2(a)\defeq  v$ such that $a = 2^v b$ for $(b,2)=1$ and $\totpos{a}$ means $a$ is totally positive.
    Then by Theorem V.1.7 in \cite{MilneCFT},  %THESIS\ref{thm:MilneExactCFT}
    since $K$ is totally real,
    \[
    h_{\mmm}=\frac{2^n(p^f-1)h}{(U:U_{\mmm,1})}.
    \]
    where $p^f\defeq \Norm_{K/\QQ}(\pp)$ for rational prime $p$ and $h\defeq h(K)=\#C_K$ is the class number of $K$ and $h_\mmm\defeq [\nRCF{\mmm}_K:K]$.

    Observe $U_{\mmm,1} \subseteq U_T$ since $\mmm$ is narrow. Then since $U_T=U^2$,
    \[(U:U_{\mmm,1})=(U:U_T)(U_T:U_{\mmm,1})=2^n(U^2:U_{\mmm,1})\]
    \[
    \implies \quad h_{\mmm}=\frac{(p^f-1)h}{(U^2:U_{\mmm,1})}.
    \]

    Consider the injection
    \[
    \frac{U^2}{U_{\mmm,1}} \hookrightarrow \left( \frac{O_K}{\pp} \right)^\times
    \]
    coming from the exact sequence and canonical isomorphism in Theorem  V.1.7 in \cite{MilneCFT}. 
    
    %THESIS\ref{thm:MilneExactCFT}.
    The image is contained in $\left(  \left( \frac{O_K}{\pp} \right)^\times\right)^2 $ so $(U^2:U_{\mmm,1}) | \frac{p^f-1}{2}$. %Then $(U_T:U_{\mm,1}) | \frac{p^f-1}{2}$. 
    Therefore $h_{\mmm}$ is even.

    %Next, we show the 2-part of the ray class group over $K$ of conductor $\mm$ is cyclic. Let $L_\pp$ denote the maximal 2-extension of the ray class field over $K$ of conductor $\mm$ where $\pp$ is a prime in $K$. Suppose $L_\pp/K$ is not cyclic and consider the inertia group at $\pp$, denoted $E$. Note that $L_\pp/K$ is tamely ramified at $\pp$ since this is a 2-extension and 2 is prime to $\pp$. This implies $E$ is cyclic so $E$ is a proper subgroup of $G=\Gal(L_\pp/K)$. Then the fixed field of $E$, denoted $L_E$, is a nontrivial even extension of $K$ in which $\pp$ is unramified, but this implies $L_E$ is contained in the ray class field over $K$ of conductor $\mm_\infty$, which is odd by Lemma \ref{RCFh} since $h(K)$ is odd so this is a contradiction.
    
    Now we show the 2-part of the ray class group over $K$ of conductor $\mmm$ is cyclic. %\switch{We use that $K$ is Galois, totally real with $U_T=U^2$ and that $K$ has odd class number $h(K)$.}
    Let $L_\pp$ denote the maximal 2-extension of the narrow ray class field over $K$ of conductor $\pp \infty$ where $\pp$ is an odd prime in $K$.
    
    Let $E$ denote the inertia group for $\pp$ relative to the extension $L_\pp/K$ and let $L_E$ denote the fixed field of $E$. 
    Note that $E$ is cyclic by Corollary 7.59 in \cite{MilneANT} because  {$L_\pp/K$ is tamely ramified at $\pp$ since $\pp$ is co-prime to $[L_\pp:K]$.}

    By Lemma \ref{RCFh}\switch{RCFh}, since $U_T=U^2$ and $K$ has odd class number, the narrow Hilbert class group has odd degree over $K$. Therefore there is no non-trivial even extension of $K$ in which all finite primes are unramified.
    Note that all finite primes of $K$ are unramified in $L_E$ and that $[L_E:K]\mid[L_\pp:K]$. Then since $[L_\pp:K]$ is a power of 2, $[L_E:K]=1$ so $E=\Gal(L_\pp/\QQ)$.

    \end{proof}

    %3.0: [defn:PrimeUnitMap] Define \gls{vK}, the map from the unramified primes of $K$ to $U/U^2$ which will induce for ``nice" $K$ a canonical isomorphism between $U/U^2$ and it's dual, a quotient of the narrow ray class group of conductor 4. 
    
    %3.1: [PrimeUnitWD] Proof that the map \gls{vK} induces a well-defined homomorphism from the narrow ray class group of conductor 4 to $U/U^2$.

%\begin{lem}[lemma:narrow2odd]\label{lemma:narrow2odd}
%Let $K$ be a totally real number field with odd class number such that $U_T=U^2$ and 2 is inert in $K/\QQ$. Then the 2-part of the narrow ray class group over $K$ of conductor 2 is trivial. 
%\end{lem}

%\begin{proof}
%Applying Theorem \ref{thm:MilneExactCFT} with  $\mm$ the narrow modulus with finite part equal to 2, the order of the narrow ray class group of conductor 2 is given by
%\[
%h_\mm = \frac{2^n h}{[U:U_{\mm,1}]}
%\]
%since setting $\mm_0=2$, a prime in $K$ by assumption,
%\[
%\Norm(\mm_0) \prod_{\pp|\mm_0} \left(1 - \frac{1}{\Norm(\pp)}\right)  = 1.
%\]
%We know $2^n | (U:U_{\mm,1})$ since $U_{\mm,1}\subseteq U_T = U^2 \subseteq U$ and $U\ism (\ZZ/2)^n$ by Dirichlet's Unit Theorem. Therefore $h_\mm | h$ so if $h$ is odd, then $h_{\mm}$ is also odd.
%\end{proof}

%%%% MISTAKE: 2 inert means Norm(2)=2^n, not 2.

\begin{lem}\label{lem:KConrad:Eisenstein}\switch{lem:KConrad:Eisenstein}
Let $L/K$ be an extension of number fields. Then $\pp$ a prime of $K$ is totally ramified in $L/K$ if and only if $L=K(\gamma)$ for some $\gamma\in \OL$ with minimal polynomial over $K$ Eisenstein at $\pp$.
\end{lem}

\switchpf{
\begin{proof}
%Localize the extension $L/K$ at the prime $\pp$ of $K$. Apply the Proposition on page 66 in Chapter III Section3 of \cite{Koblitz}. Then apply Hensel's Lemma.
See Propositions 3.10.3 and 3.10.4 in \cite{HelmutKoch}.
\end{proof}}

%%%%%%%%%%%%%%%%%

Let $\PP_K^{2\ell}$ denote the set of primes of $K$ that are co-prime to $2\ell$ where $\ell$ is the conductor of $K$.
Note that $\pp\in\PP_K^{2\ell}$ implies $\pp$ is unramified in $K/\QQ$.
%Therefore $\PP_K^{2l}=\PP_{K/\QQ}^2$, the set of primes of $K$ co-prime to 2 and unramified in $K/\QQ$.

\begin{prop}[The Treasure Map]\label{treasuremap}\switch{treasuremap}
Let $K$ be a totally real number field, Galois %\switch{(cyclic?)} 
over $\QQ$ such that the class number  $h(K)$ is odd, every totally positive unit is a square unit, and 2 is inert in $K/\QQ$. 
Define
\begin{align*}
\mathbf{x}_K:  \PP_K^{2\ell}  & \to {U/U^2} \\
  \pp &\mapsto u_\pp 
\end{align*}
such that $u_\pp\alpha_\pp$ is a square (element) in $\nRCF{\pp}_K$ where $\alpha_\pp\in\OK$ is a totally positive generator of $\pp^{h(K)}$. The map $\mathbf{x}_K$ is well-defined.
\end{prop}
\switchpf{
    \begin{proof} 
Let $\pp\in\PP_K^{2\ell}$.
By Lemma \ref{Uquadsub}\switch{Uquadsub}, %K tot real UT=U^2, Galois, odd class number,
there is a unique quadratic extension $L/K$ such that $L\subseteq\nRCF{\pp}_K$. 
%the narrow ray class field over $\QQ$ of conductor $\pp$. 
By Lemma \ref{lem:KConrad:Eisenstein}\switch{lem:KConrad:Eisenstein}, 
\[
L=K(\gamma)
\]
for some $\gamma\in \OL$ with minimal polynomial over $K$ Eisenstein at $\pp$. 

%We may assume $\Trace_{L/K}(\gamma)=0$; otherwise taking $\gamma'\defeq 2\gamma - \Trace_{L/K}(\gamma)$, the minimal polynomial of $\gamma'$ over $K$ is also Eisenstein at $\pp$ and $\Trace_{L/K}(\gamma')=0$. 

Let $f_\gamma(x)=x^2 + c_1 x + c_0\in K[x]$ be the minimal polynomial of $\gamma$ over $K$. We may assume $c_1=0$ since otherwise, we could take $\gamma'\defeq  2\gamma + c_1\in L$ and the minimal polynomial of $\gamma'$ is $x^2 - (c_1^2-4c_0)$, which is also Eisenstein at $\pp$ since $\pp$ is odd.
Now write the minimal polynomial of $\gamma$ over $K$ as
\[
f_\gamma(x) = x^2 - c
\]
for some $c\in K$ where $f_\gamma(x)$ is Eistenstien at $\pp$.

%We may assume $\gamma\in\mathcal{O}_L$; otherwise $c= \frac{a}{b}$ for some $a,b\in\OK$ and $\ord_\pp(c)=1$ implies $\ord_\pp(a)=1$ and $\ord_\pp(b)=0$. Then $\gamma' \defeq  b\gamma$ has minimal polynomial $x^2 -ab$, which is Eisenstein at $\pp$ and $ab\in\OK$ implies $\gamma'\in\mathcal{O}_L$.

%Note that $f_\gamma(x)$ can not be Eisenstein at any other primes of $K$ by Lemma \ref{lem:KConrad:Eisenstein}\switch{lem:KConrad:Eisenstein} because $L\subseteq \nRCF{\pp}_K$ implies $\pp$ is the only finite prime that ramifies in $L/K$.

If $c\OK$ had a prime factor $\qq\neq \pp$ with odd multiplicity then $L$ would be ramified at $\qq$, a contradiction since $L\subseteq \nRCF{\pp}_K $. Therefore
\[
c\OK = \pp I^2
\]
for some ideal $I\subseteq \OK$ co-prime to $\pp$.

Let $b\in\OK$ be a generator of $I^{h(K)}$. 
% b^2 is well-defined up to mult. by a square unit
Raising $c\OK = \pp I^2$ to the power of $h(K)$ gives
\[
\frac{c^{h(K)}}{b^2} \OK = \pp^{h(K)}.
\]
Let $d = {c^{h(K)}}/{b^2}$. Then $d\OK = \pp^{h(K)}$.
Since $\pp^{h(K)}$ is principal, $d\in\OK$. Furthermore, since $h(K)$ is odd, 
\[
d = c \left( \frac{c^{\frac{h(K)-1}{2}}}{b} \right)^2
%\implies
\]
and so $K(\sqrt{d}) = K(\sqrt{c})$.
We have proven that $L=K(\sqrt{d})$ for some generator $d\in\OK$ of $\pp^{h(K)}$. 

By Lemma \ref{RCFh}\switch{RCFh}, since $U_T=U^2$, there exists a totally positive generator $\alpha\in\OK$ of $\pp^{h(K)}$.
There exists some $u\in \OK^\times$ such that
\[d = u \alpha\] 
because all signatures are represented by units since $U_T=U^2$ in $K$. 
Since $L=K(\sqrt{d})\subseteq\nRCF{\pp}_K$, $d=u\alpha$ is a square (element) in $\nRCF{\pp}_K$. 

If $u,v\in\OK^\times$ such that $u\alpha$ and $v\alpha$ are both squares in $\nRCF{\pp}_K$, then
$K(\sqrt{u\alpha})$ and $K(\sqrt{v\alpha})$ are both contained in $\nRCF{\pp}_K$. By uniqueness of $L$, $K(\sqrt{u\alpha})=K(\sqrt{v\alpha})$.
Write $\sqrt{v\alpha} = x + y \sqrt{u\alpha}$. Squaring both sides yields $xy=0$, which implies $x=0$ since $y\neq 0$. Thus $v\alpha = y^2 u \alpha$. Since ${y^2}$ is totally positive, $u$ and $v$ have the same signature. Since $U_T=U^2$, this implies $u\equiv v$ in $U/U^2$.
    \end{proof}
}%end switchpf
%OLd PROOFS are in research journal Chap3 .tex file,  as a \switch{} comment.

%%%%%%%%%%%%%%

%%%%%%%%%%%%%%%%%%%%%%%

\section{Density Formulas Modulo 4}\label{sec:pmdensityformulas:orange}%[thm:pmdensityformulas]%[spindep]%[pmstarwelldefined]%[StarMG]%[defn:NMG]

%\begin{center}
%    \switch{sec:pmdensityformulas:orange}
%\end{center}

%%%%%%%%%%%%%%

%\section{Equidistribution and the Starlight Invariants}\label{sec:equidistribution}

In this section, we prove a more refined version of the main results of \cite{APSR} that handles restriction modulo $4\ZZ$.

Lemma 11.1 in \cite{FIMR} states that the product $\spin(\pp,\sigma)\spin(\pp,\sigma\inv)$ is a product of Hilbert symbols at places dividing 2. We restate this more explicitly in Lemma \ref{spindep}\switch{spindep}.
For a place $v$ of $K$, let $K_{(v)}$ denote the completion of $K$ at $v$. 
For $a,b\in K$ co-prime to $v$, the Hilbert symbol is defined such that $(a,b)_v\defeq 1$ if the equation $ax^2+by^2=z^2$ has a solution $x,y,z\in K_{(v)}$ where at least one of $x$, $y$, or $z$ is nonzero and $(a,b)_v\defeq -1$ otherwise. 

\begin{lem}\cite{FIMR}\label{spindep}\switch{spindep}
Let $K$ be a cyclic number field of degree $\geq 3$ with odd class number $h(K)$ such that every totally positive unit is a square.
Let $\alpha$ be a totally positive generator of the odd prime ideal $\pp\subseteq\OK$. Then
\[
\spin(\pp,\sigma)\spin(\pp,\sigma\inv) = \prod_{v\mid2}(\alpha,\alpha^\sigma)_v.
\]
In particular, if $\alpha\equiv 1 \bmod 4$ then $\prod_{v\mid2}(\alpha,\alpha^\sigma)_v=1$.
\end{lem}

\switchpf{
\begin{proof}
See Lemma 11.1 in \cite{FIMR} or use the fact that Hilbert symbols satisfy $\prod_{v}(\alpha,\alpha^\sigma)_v=1$.
\end{proof}
}

%\begin{rmk}\label{rmk:2inert:spindep}\switch{rmk:2inert:spindep}
%When 2 is inert in $K/\QQ$, this becomes
%\[
%\spin(\pp,\sigma)\spin(\pp,\sigma\inv) = (\alpha,\alpha^\sigma)_2.
%\]
%\end{rmk}
%%%%%%%%%%%%%%

For a number field $K$ that is Galois over $\QQ$ with ring of integers $\OK$, %such that 2 is inert in $K/\QQ$, 
define
\[
\MK\defeq  (\OK/4\OK)^\times / \left( (\OK/4\OK)^\times \right)^2.
\]
Let $\mathbf{M}_{4,G}$ denote the set of $\Gal(K/\QQ)$-orbits of $\MK$.  The following definition is motivated by Lemma \ref{spindep}\switch{spindep}.

\begin{lem}\cite{APSR}
\label{StarMG}\switch{StarMG}
Let $K$ be a number field, Galois over $\QQ$ with abelian Galois group. %such that 2 is inert in $K/\QQ$.
Let $\alpha\in\OK$ denote a
representative of $[\alpha]\in$ $\mathbb{M}_{4,G}$. 
Define the map
\begin{align*}
\text{$\star$}:  \text{ $\mathbb{M}_{4,G}$ }& \to \{\pm 1\} \\
  [\alpha] &\mapsto \left\{ 
  \begin{array}{l l}
    1 & \quad \text{if }  %(\alpha,\alpha^\sigma)_2 = 1
    \prod_{v2}(\alpha,\alpha^\sigma)_v = 1 \text{ for all non-trivial } \sigma\in \Gal(K/\QQ),  \\
    -1 & \quad \text{otherwise.}\\
  \end{array} \right.
\end{align*}
Then $\star$ is a well-defined map.
\end{lem}

%Note that by Lemma \ref{spindep}, a rational prime $p$ satisfies the spin relation 
%\[
%\spin(\pp,\sigma) = \spin(\pp,\sigma\inv) \quad \text{for all } \sigma\neq 1 \in \Gal(K/\QQ)
%\]
%where $\pp$ is a prime of $K$ above $p$ exactly when $\star\circ\mathbf{r}(p)=1$ where $\mathbf{r}$ is as defined in Definition \ref{defn:r}.

\switchpf{
\begin{proof}
This is Theorem 5.1\switch{(StarMG)} in \cite{APSR}.
\end{proof}
}

%%%%%%%%%%%%%%

Let $\PP_K^{2}$ denote the set of primes of $K$ that are co-prime to $2$. %Recall the following maps from \cite{APSR} which were shown to be well-defined following Definition 4.4. %in \cite{APSR}. 

\begin{defn}\cite{APSR}\label{defn:r} \switch{[defn:r]}
Let $q\geq 4$ be a power of $2$. 
Let $K$ be a totally real number field such that every totally positive unit is a square unit. Assume $K$ is Galois over $\QQ$ with abelian Galois group and that $K$ has odd class number, $h(K)$.
\begin{enumerate}
    \item  
    Define the map
\begin{align*}
\mathbf{r}_0:  
\PP_K^{2}&\to  \mathbf{M}_q \\
  \pp &\mapsto \alpha,
\end{align*}
where $\alpha\in\OK$ is a totally positive generator for the principal ideal $\pp^{\text{$h(K)$}}$. 
%\PP_K^{2\ell} was called S' in the journal
\item %Let $\PP_\QQ^{2}$ denote the set of odd rational primes.  
Define the map
\begin{align*}
\mathbf{r}:  
\PP_\QQ^{2}&\to \text{ $\mathbb{M}_{q,G}$ } \\
  p &\mapsto [\mathbf{r}_0(\pp)],
\end{align*}
where $\pp$ is any prime in $K$ above $p$. 
Here $[\alpha]$ denotes the $\Gal(K/\QQ)$-orbit of $\alpha\in \mathbf{M}_4$ considered in $\mathbb{M}_{q,G}$. 
\end{enumerate}
\end{defn}
\switchpf{
\begin{proof}
To see that $\mathbf{r}_0$ is well-defined, apply Lemma \ref{RCFh}\switch{RCFh} which gives the existence of a totally positive generator  $\alpha$ for the principal ideal $\pp^{h(K)}$. Since $U_T=U^2$ and squares are trivial in $\mathbf{M}_q$, $\mathbf{r}_0$ is well-defined.
Next observe that $\mathbf{r}_0(\pp^\sigma) =  \mathbf{r}_0(\pp)^\sigma$ for all $\pp\in \PP_K^2$ and all $\sigma\in \Gal(K/\QQ)$. Then since $\mathbb{M}_{4,G}$ is the set of Galois orbits of $\mathbf{M}_4$, 
$
[\mathbf{r}_0(\pp^\sigma)] = [\mathbf{r}_0(\pp)] 
$
in $\mathbb{M}_{4,G}$ and so $\mathbf{r}$ is well-defined.
\end{proof}
}

\begin{defn}\label{defn:density} Let $S$ be a set of primes and let $R\subseteq S$. If the limit exists, then the \textit{restricted density} of $R$ (restricted to $S$) is defined as
\[
d(R|S) = d\left( \frac{R}{S} \right)\defeq  \lim_{N\to \infty} \frac{\# R_N}{\#S_N}\]
where $S_N$ and $R_N$ denote the set of primes in $S$ and $R$ respectively of absolute norm less than $N\in \ZZ_+$. 
\end{defn}

We now restate Lemma 4.3\switch{(equidistribution APSR)} from \cite{APSR}. %that handles the densities restricted to primes of a fixed congruence class modulo $4\ZZ$. 
Let $\PP_K^{2\ell}$ denote the set of primes of $K$ that are co-prime to $2\ell$ where $\ell$ is the conductor of $K$.
 By Lemma 3.5 \switch{(homMK:surj)} in \cite{APSR}, the map $\mathbf{r}_0$ from Definition \ref{defn:r}\switch{defn:r} induces a canonical surjective homomorphism
\[\varphi: \nRCG{4}_K\twoheadrightarrow \MK.\]

\begin{lem}\cite{APSR}\label{equidistribution}\switch{[equidistribution]} Let $K$ be a cyclic totally real number field with odd class number such that every totally positive unit is a square. Assume $2$ is inert in $K/\QQ$.  Assume $n\defeq [K:\QQ]$ is an odd prime. Let $\ell$ denote the conductor of $K$. 

\begin{enumerate}
    \item\label{equidistribution:allp}  For any $\alpha\in \mathbf{M}_4$, the density of $\pp\in \PP_K^{2\ell}$ such that $\varphi(\pp)=\alpha$ is ${1}/{2^n}$. That is, 
    \[
    d(\mathbf{r}_0\inv(\alpha)| \PP_K^{2\ell}) = \frac{1}{\#\mathbf{M}_4} = \frac{1}{2^n}.
    \]
    \item\label{equidistribution:spcomp}  Furthermore, the density does not change when we restrict to primes of $K$ that split completely in $K/\QQ$. That is, %if $S'$ is the set of primes $\pp\in\PP_K^{2\ell}$ with \gls{inertia degree} $f_{K/\QQ}(\pp)=1$ in $K/\QQ$, then
    \[
    d(\mathbf{r}_0\inv(\alpha) \cap S' | S') = \frac{1}{\#\mathbf{M}_4} = {1}/{2^n}.
    \]
\end{enumerate}
\end{lem}

\begin{proof}
For part (1), see Lemma 4.3\switch{equidistribution} in \cite{APSR}.
We now prove part (2). 
By Lemma 3.5\switch{(homMK:surj)} in \cite{APSR}, the map $\mathbf{r}_0: \PP_K^{2\ell} \to \MK$ from Definition \ref{defn:r} induces a surjective canonical group homomorphism
\[
\varphi: \nRCG{4}_K \twoheadrightarrow \MK
\]
which commutes with the Galois action from $\Gal(K/\QQ)$. 
Define $H\defeq  \Art(\ker(\varphi))$ so that the following diagram of exact sequences commutes and commutes with the Galois action from $\Gal(K/\QQ)$. 
\begin{center}
\begin{tikzpicture}
\diagram (m)
{ 1 & \ker(\varphi) & \nRCG{4}_K & \mathbf{M}_4 & 1  \\
  1 & H &  \Gal(\nRCF{4}_K/K) &  \mathbf{M}_4 & 1 \\};
\path [->] 
	   %top
           (m-1-1) edge node [above] {} (m-1-2)
           (m-1-2) edge node [above] {} (m-1-3)
           (m-1-3) edge node [above] {$\varphi$} (m-1-4)
           (m-1-4) edge node [above] {} (m-1-5)
           %bottom
           (m-2-1) edge node [above] {} (m-2-2)
           (m-2-2) edge node [above] {} (m-2-3)
           (m-2-3) edge node [above] {} (m-2-4)
           (m-2-4) edge node [above] {} (m-2-5)
           ;
\path [->] 
           %vertical
           (m-1-2) edge node [right] {$\Art$} (m-2-2)
           (m-1-3) edge node [right] {$\Art$} (m-2-3)
           (m-1-4) edge node [right] {\text{identity}} (m-2-4)
          ;
\end{tikzpicture}
\end{center}
Define $L$ to be the fixed field of $H$. Then $\varphi$ induces a canonical isomorphism 
\[
\Gal(L/K)\ism \MK.
\]

For a Galois extension $F/E$ of conductor dividing $\mmm$, let
\begin{align*}
\preArt{F|E}{E}(\tau)&\defeq \{ p\in\PP_E^{\mmm}: \Art_{F|E}(p) = \tau\},\\
\preArt{F|E}{F}(\tau)&\defeq \{ \pp\in\PP_F^{\mmm}: \pp \text{ lies above }p\in\preArt{F|E}{E} \},
%\Art_{F|E}(\pp) = \tau\},
\end{align*}
where $\tau\in\Gal(F/E)$.

Let $\alpha\in \MK$.
Let $\sigma\in\Gal(L/K)$ corresponding to $\alpha\in\MK$. Then $\preArt{L|K}{K}(\sigma) = \mathbf{r}_0\inv(\alpha)$ is the set of primes of $K$ that map to $\alpha$ via $\mathbf{r}_0$ and $\preArt{K|\QQ}{K}(1) =S'$ is exactly the set of primes of $K$ that split completely from $K/\QQ$.
We want to show
\[
d\defeq d( \preArt{L|K}{K}(\sigma) \cap \preArt{K|\QQ}{K}(1) | \preArt{K|\QQ}{K}(1)) = \frac{1}{2^n}.
\]

Since the diagram above commutes with the action from $\Gal(K/\QQ)$, $L$ is Galois over $\QQ$ with $\Gal(L/K)\trianglelefteq \Gal(L/\QQ)$. Considering $\sigma\in\Gal(L/\QQ)$, since $\sigma$ fixes $K$, 
$
\preArt{L|K}{K}(\sigma) \subseteq \preArt{K|\QQ}{K}(1).
$
Therefore 
$
d=d( \preArt{L|K}{K}(\sigma)| \preArt{K|\QQ}{K}(1)).
$

For all $N\to\infty$, there is a surjective map
\[
\preArt{K/\QQ}{K}(1)|_N \twoheadrightarrow \preArt{K/\QQ}{\QQ}(1)|_N
\]
with index $n=[K:\QQ]$, where $|_N$ denotes the restriction to primes of absolute norm less than $N$. There is also a surjective map
\[
 \preArt{L|K}{K}(\sigma)|_N \twoheadrightarrow \preArt{L|\QQ}{\QQ}(\sigma)|_N
\]
of index $\#\Stab(\sigma)$, where $\Stab(\sigma)$ is the Stabilizer of $\sigma\in\Gal(L/K)$ under the action of $\Gal(K/\QQ)$. Therefore
\[
\# \preArt{L/\QQ}{\QQ}(\sigma)|_N =\frac{\# \preArt{L/K}{K}(\sigma)|_N }{\#\Stab(\sigma)}
\quad \text{and} \quad 
\#\preArt{K/\QQ}{\QQ}(1)|_N = \frac{\#\preArt{K/\QQ}{K}(1)|_N}{\#\Gal(K/\QQ)}. 
\]
Therefore
\[
d = d( \preArt{L/\QQ}{\QQ}(\sigma)|\preArt{K/\QQ}{\QQ}(1) ) \frac{\#\Stab(\sigma)}{\#\Gal(K/\QQ)}.
\]
Therefore, by Cebotarev's Theorem (see Theorem VII.13.4 in \cite{NeukirchANT} for Dirichlet density or Theorem 4 in \cite{Serre} for natural density),
\begin{align*}
d &= 
\frac{
\left(\frac{\#\Stab(\sigma)\#\left< \sigma \right>}{\#\Gal(K/\QQ)\#\Gal(L/\QQ)}\right)
}{
\left(\frac{1}{\#\Gal(K/\QQ)}\right)
}\\
&=
\frac{\#\Stab(\sigma)\#\left< \sigma \right>}{\#\Gal(K/\QQ)\#\Gal(L/K)},
\end{align*}
where $\left< \sigma \right>$ denotes the conjugacy class of $\sigma$ under the action of $\Gal(K/\QQ)$. By the Orbit-Stabilizer Theorem, $\#\Stab(\sigma)\#\left< \sigma \right> = \#\Gal(K/\QQ)$. Therefore
\[
d = \frac{1}{\#\Gal(L/K)} = \frac{1}{2^n},
\]
since $\Gal(L/K) \ism \MK$ and $\#\MK=2^n$ by Lemma 3.3\switch{(Mundy)} in \cite{APSR}.

\end{proof}

By Lemma \ref{equidistribution}\switch{equidistribution}, for any $\alpha \in\mathbf{M}_4$, the density of primes of $K$ that map to $\alpha$ via $\mathbf{r}_0$ is $\frac{1}{\#\mathbf{M}_4}= \frac{1}{2^n}$. That is, the primes of $K$ are equidistributed in $\mathbf{M}_4$ via the map $\mathbf{r}_0$. We apply this Lemma in the following definition and then we state a more refined version of this Lemma. 

\begin{defn}\label{defn:NMG}\switch{defn:NMG}
Let $K$ be a cyclic totally real number field with odd class number such that every totally positive unit is a square. For simplicity, we assume 2 is inert in $K/\QQ$.  Assume $n\defeq [K:\QQ]$ is an odd prime. Let $\ell$ denote the conductor of $K$.

Let $\alpha\in\MK$. Let $\pp$ be a prime in $\PP_K^{2\ell}$ such that $\mathbf{r}_0(\pp)=\alpha$. The map 
\begin{align*}
\mathbf{N}: \MK &\to (\ZZ/4)^\times \\
  \alpha &\mapsto \Norm_{K/\QQ}(\pp) \bmod 4\ZZ
\end{align*}
is well-defined and $\mathbf{N}(\pp)=\mathbf{N}(\pp^\sigma)$ for all $\sigma\in\Gal(K/\QQ)$. 
%\begin{align*}
%\mathbf{N}: \MKG &\to (\ZZ/4)^\times \\
% X &\mapsto p \bmod 4\ZZ
%\end{align*}
%is well-defined.
\end{defn}
\switchpf{
\begin{proof}
 By  Lemma \ref{equidistribution}\switch{equidistribution}, for any $\alpha\in \MK$, there exists a prime $p\in\PP_\QQ^{2\ell}$ such that $\mathbf{r}(p)=\alpha$. 
 
 Let $p,q\in\PP_\QQ^{2\ell}$ such that $\mathbf{r}(p)=\mathbf{r}(q)$. We will show that $p\equiv q \bmod 4\ZZ$. Let $\pp,\qq\in\PP_K^{2\ell}$ be primes above $p$ and $q$ respectively with totally positive generators $\alpha,\beta\in\OK$ respectively. Since $\mathbf{r}(p)=\mathbf{r}(q)$, 
 \[
 \mathbf{r}_0(\pp)^\tau=\mathbf{r}_0(\qq)
 \]
 for some $\tau\in\Gal(K/\QQ)$. Therefore 
 \[
 \alpha^\tau \equiv \beta y^2 \bmod 4\OK
 \]
 for some $y\in\OK$ %(y is co-prime to 2.)
 and so
 $\alpha^{\tau\sigma} \equiv \beta^\sigma (y^\sigma)^2 \bmod 4\OK$  for all  $\sigma\in\Gal(K/\QQ)$, which implies
 \[
\Norm_{K/\QQ}(\alpha) \equiv \Norm_{K/\QQ}(\beta) \bmod 4\OK.
\]

 Write $\Norm_{K/\QQ}(\alpha)- \Norm_{K/\QQ}(\beta) = 4 \gamma$ where $\gamma\in\OK$. Then since the norms are in $\ZZ$, we know $4\gamma\in\ZZ$ so $\gamma\in \OK \cap \QQ=\ZZ$. Therefore
 \[
 \Norm_{K/\QQ}(\alpha) \equiv \Norm_{K/\QQ}(\beta) \bmod 4\ZZ.
 \]
 If $p$ and $q$ split completely in $K/\QQ$ we are done since $p=\Norm_{K/\QQ}(\alpha)$ and $q=\Norm_{K/\QQ}(\beta)$. Otherwise since $K/\QQ$ is cyclic, $p$ or $q$ is inert. 
 
 If $p$ is inert, $p^n=\Norm_{K/\QQ}(\alpha)$ where $n=[K:\QQ]$. Recalling that $n$ is odd and $p\neq2$,
 \[
 \left( \frac{p^{n-1}}{4}\right) =1,
 \]
 where the left hand side is the quadratic residue symbol in $\ZZ$. Therefore $p^{n-1} \equiv 1 \bmod 4$ so $p^n \equiv p \bmod 4$. Thus no matter how $p$ factors in $K/\QQ$, 
 \[
 \Norm_{K/\QQ}(\alpha)\equiv p \bmod 4\ZZ
 \]
 and the analogous statement is true for $q$.
 Therefore
 \[
 p \equiv q \bmod 4\ZZ.
 \]
 \end{proof}}

We now state an extended version of Lemma \ref{equidistribution}\switch{equidistribution} that handles the densities restricted to primes of a fixed congruence class modulo $4\ZZ$. Recall that $\PP_K^{2\ell}$ denotes the set of primes of $K$ that are co-prime to $2\ell$ where $\ell$ is the conductor of $K$.
 Recall that by Lemma 3.5\switch{(homMK:surj)} in \cite{APSR}, the map $\mathbf{r}_0$ from Definition \ref{defn:r}\switch{defn:r} induces a canonical surjective homomorphism,
\[\varphi: \nRCG{4}_K\twoheadrightarrow \MK.\]

\begin{lem}\label{equidistribution:extpm}\switch{[equidistribution:extpm]} Let $K$ be a cyclic totally real number field of odd prime degree $n=[K:\QQ]$ with odd class number such that every totally positive unit is a square. Assume $2$ is inert in $K/\QQ$.  Let $\ell$ denote the conductor of $K$.

    For a fixed sign $\pm$, let $S_{\pm}'$ denote the set of primes of $K$ laying above some $p\in S$ such that $p\equiv \pm 1 \bmod 4\ZZ$.
    For any $\alpha\in \mathbf{M}_4$, the density of $\pp\in S_{\pm}'$ such that $\varphi(\pp)=\alpha$ is given by
    \[
    d(\mathbf{r}_{0}\inv(\alpha) \cap S_{\pm}'|S_{\pm}') = \left\{
    \begin{array}{ll}
        \frac{1}{2^{n-1}} &  \text{if } \mathbf{N}(\alpha) = \pm 1 \bmod 4\\
        0 & \text{otherwise.}
    \end{array}
    \right.
    \]
%\end{enumerate}
\end{lem}

\begin{proof}
%%\todo{To Do: Re-write this proof. This proof does not show that the limit defining the density in question exists.}
By Lemma 3.5\switch{(homMK:surj)} in \cite{APSR}, the map $\mathbf{r}_0: \PP_K^{2\ell} \to \MK$ from Definition \ref{defn:r}\switch{defn:r} induces a surjective canonical group homomorphism
\[
\varphi: \nRCG{4}_K \surj \MK,
\]
which commutes with the Galois action from $\Gal(K/\QQ)$. 
Define $H\defeq  \Art(\ker(\varphi))$ so that the following diagram of exact sequences commutes and commutes with the Galois action from $\Gal(K/\QQ)$. 

\begin{center}
\begin{tikzpicture}
\diagram (m)
{ 1 & \ker(\varphi) & \nRCG{4}_K & \mathbf{M}_4 & 1  \\
  1 & H &  \Gal(\nRCF{4}_K/K) &  \mathbf{M}_4 & 1 \\};
\path [->] 
	   %top
           (m-1-1) edge node [above] {} (m-1-2)
           (m-1-2) edge node [above] {} (m-1-3)
           (m-1-3) edge node [above] {$\varphi$} (m-1-4)
           (m-1-4) edge node [above] {} (m-1-5)
           %bottom
           (m-2-1) edge node [above] {} (m-2-2)
           (m-2-2) edge node [above] {} (m-2-3)
           (m-2-3) edge node [above] {} (m-2-4)
           (m-2-4) edge node [above] {} (m-2-5)
           ;
\path [->] 
           %vertical
           (m-1-2) edge node [right] {$\Art$} (m-2-2)
           (m-1-3) edge node [right] {$\Art$} (m-2-3)
           (m-1-4) edge node [right] {\text{identity}} (m-2-4)
          ;
\end{tikzpicture}
\end{center}
Define $L$ to be the fixed field of $H$ 
in $\Gal(\nRCF{4}_K/K)$. Then $\varphi$ induces a canonical isomorphism 
\[
\Gal(L/K)\ism \MK,
\]
which commutes the with action from $\Gal(K/\QQ)$.

By Theorem 1.7 in \cite{MilneCFT}, since $U_T=U^2$ and 2 is inert in $K/\QQ$, 
\[[\nRCF{4}_K:K] \mid h2^n(2^n-1)\]
where $h$ is the class number of $K$. We know $[L:K]=2^n$ since $\Gal(L/K)\ism \mathbf{M}_4$ and $\#\MK=2^n$ by Lemma 3.3\switch{(Mundy)} in \cite{APSR}. Therefore $[\nRCF{4}_K:L]$ is odd. Letting $K_4$ denote the composite of $K$ and $\QQ(\zeta_4)$, this implies that $K_4\subseteq L$ so we have the following tower of number fields of the following degrees.

\begin{center}
\begin{tikzpicture}
\diagram (m)
{ L \\
   K_4 \\
   K \\
   \QQ \\};
\path [-] 
 	(m-1-1) edge node [right]{$2^{n-1}$} (m-2-1)
	(m-2-1) edge node [right]{$2$} (m-3-1)
	(m-3-1) edge node [right]{$n$} (m-4-1)
	;
\end{tikzpicture}
\end{center}

For a Galois extension $F/E$ of conductor dividing $\mmm$, let
\begin{align*}
\preArt{F|E}{E}(\tau)&\defeq \{ p\in\PP_E^{\mmm}: \Art_{F|E}(p) = \tau\},\\
\preArt{F|E}{F}(\tau)&\defeq \{ \pp\in\PP_F^{\mmm}: \pp \text{ lies above }p\in\preArt{F|E}{E}(\tau) \},
%\Art_{F|E}(\pp) = \tau\},
\end{align*}
where $\tau\in\Gal(F/E)$.

Let $\alpha\in\MK$ and let $\sigma\in\Gal(L/K)$ corresponding to $\alpha$ via the isomorphism induced by $\varphi$. Note that $\Gal(L/K)\trianglelefteq \Gal(L/\QQ)$. %We will consider $\sigma$ inside the larger Galois group, $\Gal(L/\QQ)$ such that $\sigma$ fixes $K$. %and we will use the notation $\sigma_0$ when we are considering $\sigma$ inside $\Gal(L/K)$.

Fix a sign $\pm$ and let $\tau_0\in\Gal(\QQ(\zeta_4)/\QQ)$ such that 
\[
\preArt{\QQ(\zeta_4)/\QQ}{\QQ}(\tau_0) = \{ p \in \PP_\QQ^{2\ell}: p \equiv \pm 1 \bmod 4\ZZ\}.
\]
Note that since $n=[K:\QQ]$ is odd and $[\QQ(\zeta_4):\QQ]=2$,
$\Gal(K_4/K) \ism \Gal(\QQ(\zeta_4)/\QQ)$ canonically. Let $\tau\in\Gal(K_4/K)$ corresponding to $\tau_0$ in $\Gal(\QQ(\zeta_4)/\QQ)$.
Note that $\Gal(K_4/K)\trianglelefteq\Gal(K_4/\QQ)$. 
%We use the notation $\tau'$ to refer to $\tau$ considered inside the larger group, $\Gal(K_4/\QQ)$ such that $\tau'$ fixes $K$.

Observe that 
\[
    \mathbf{r}_0\inv(\alpha) = \preArt{L/K}{K}(\sigma), \quad 
    S' = \preArt{K/\QQ}{K}(1), \quad \text{and} \quad
    S_\pm' = \preArt{K_4/K}{K}(\tau) \cap S'.
\]
Then the density in question is
\[
d_\pm\defeq  d\left( \frac{\mathbf{r}_0\inv(\alpha) \cap S_\pm' }{S_\pm'} \right)
= d \left( \frac{\preArt{K/\QQ}{K}(1) \cap \preArt{L/K}{K}(\sigma) \cap \preArt{K_4/K}{K}(\tau)}{\preArt{K/\QQ}{K}(1) \cap \preArt{K_4/K}{K}(\tau)} \right)
\]

Consider $\bar{\sigma}\in \Gal(K_4/K)$ taken to be the image of $\sigma$ under the natural surjection,
\[
\Gal(L/K) \twoheadrightarrow \Gal(K_4/K).
\]
Note that $\bar{\sigma} = \tau$ exactly when $\mathbf{N}(\alpha) = \pm 1 \bmod 4$. 
If $\bar{\sigma} \neq \tau$ then 
$
\mathbf{r}_0\inv(\alpha) \cap S_\pm' = \emptyset
$
so the density in question is 0. We now assume $\mathbf{N}(\alpha) = \pm 1 \bmod 4$ so that $\bar{\sigma} = \tau$. Then 
\[
\preArt{L/K}{K}(\sigma) \cap \preArt{K_4/K}{K}(\tau) = \preArt{L/K}{K}(\sigma).
\]
Therefore
\[
d_\pm = d \left( \frac{\preArt{K/\QQ}{K}(1) \cap \preArt{L/K}{K}(\sigma) }{\preArt{K/\QQ}{K}(1) \cap \preArt{K_4/K}{K}(\tau)} \right).
\]
Restricting to primes of norm over $\QQ$ less than $N$, there are surjective maps of the following indices
\[
\preArt{K/\QQ}{K}(1) \cap \preArt{L/K}{K}(\sigma) |_N \rightarrow \preArt{L/\QQ}{\QQ}(\sigma)|_N \quad \text{ with index = } \#\Stab(\sigma)
\]
and
\[
\preArt{K/\QQ}{K}(1) \cap \preArt{K_4/K}{K}(\tau) |_N \rightarrow \preArt{K_4/\QQ}{\QQ}(\tau)|_N \quad \text{ with index = } \#\Stab(\tau)
\]
where $\Stab(\sigma)$ denotes the stabalizer of $\sigma$ under the action from $\Gal(K/\QQ)$ and $\Stab(\tau)$ denotes the stabalizer of $\tau$ under the action from the same group. Let $\left<\sigma\right>$ denote the conjugacy class of $\sigma$ under the same action and define $\left<\tau\right>$ similarly.
Then by Cebotarev's Theorem (see Theorem VII.13.4 in \cite{NeukirchANT} for Dirichlet density or Theorem 4 in \cite{Serre} for natural density) and by the Orbit-Stabilizer Theorem,
\begin{align*}
d_\pm &= \frac{\#\Stab(\sigma)}{\#\Stab(\tau)}
d \left(\frac{
\preArt{L/\QQ}{\QQ}(\sigma)
}{
\preArt{K_4/\QQ}{\QQ}(\tau)
}
\right) \\
& = \frac{\#\Stab(\sigma) d(\preArt{L/\QQ}{\QQ}(\sigma)) }{\#\Stab(\tau) d(\preArt{K_4/\QQ}{\QQ}(\tau))} 
\\
& = \left(\frac{\#\Stab(\sigma)\#\left<\sigma\right>}{\#\Stab(\tau)\#\Gal(L/K)}\right)
\bigg/
 \left(\frac{\#\left<\tau\right>}{\#\Gal(K_4/K)}\right)
 \\
& = \frac{\#\Stab(\sigma)\#\left<\sigma\right>}{\#\Stab(\tau)\#\left<\tau\right>\#\Gal(L/K_4)}  \\
& = \frac{1}{\#\Gal(L/K_4)} \\
&= \frac{1}{2^{n-1}}.
\end{align*}

\end{proof}

%\newpage
\begin{defn}\label{defn:StarlightInvariantpm}\switch{defn:StarlightInvariantpm} Let $K$ be a totally real cyclic number field of odd prime degree $n$ with odd class number $h(K)$ such that every totally positive unit of $K$ is a square. Assume 2 is inert in $K$. %Let $\ell$ denote the conductor of $K$.
We define the \textit{positive Starlight invariant} of $K$ to be 
\[
\mKpos\defeq \#\{X\in \MKG: \#X=n, \star(X)=1,  \text{ and } \mathbf{N}(X)=1 \}
\]
and we define the \textit{negative Starlight invariant} of $K$ to be 
\[
\mKneg\defeq \#\{X\in \MKG: \#X=n, \star(X)=1,  \text{ and } \mathbf{N}(X)=-1 \}.
\]
\end{defn}

\begin{rmk}
By Lemma 5.2 \switch{pm1} and Lemma 3.3 \switch{Mundy} in \cite{APSR}, 
\[
\#\ker(\star) = \mKpos + \mKneg + 1.
\]
\end{rmk}

%%%%%%%%%%%%%%

%\subsection{Density formulas Restricted Modulo 4}

Define $S_+$ to be the set of odd rational primes congruent to $1 \mod 4$ that split completely in $K/\QQ$. Similarly, define $S_-$ to be the set of odd rational primes congruent to $-1 \mod 4$ that split completely in $K/\QQ$. 
For a fixed sign $\pm$, define
\[
R_\pm \defeq \{p\in S_\pm: \spin(\pp,\sigma) = \spin(\pp,\sigma\inv)  \text{ for all } \sigma\neq 1 \in \Gal(K/\QQ) \}. 
\]
%\begin{align*}
    %& R_+ \defeq \{p\in S_+: \spin(\pp,\sigma) = \spin(\pp,\sigma\inv)  \text{ for all } \sigma\neq 1 \in \Gal(K/\QQ) \} \quad \text{and}\\
    %& R_- \defeq \{p\in S_-: \spin(\pp,\sigma) = \spin(\pp,\sigma\inv)  \text{ for all } \sigma \neq 1 \in \Gal(K/\QQ) \},
%\end{align*}

\begin{thm}[$\pm$ Density Formulas]\label{thm:pmdensityformulas}\switch{thm:pmdensityformulas}
Let $K$ be a cyclic totally real number field of odd prime degree $n$ with odd class number such that 2 is inert in $K/\QQ$ and every totally positive unit is a square unit. Then
\[
d_+\defeq d(R_+|S_+) = \frac{1+\mKpos n}{2^{n-1}} \quad \text{and} \quad
d_-\defeq d(R_-|S_-) = \frac{\mKneg n}{2^{n-1}}.
\]
\end{thm}

\switchpf{
\begin{proof} 

Recall the map $\mathbf{r}$ from Definition \ref{defn:r}\switch{defn:r}.
 By Lemma \ref{spindep}\switch{spindep}, for each fixed sign $\pm$,
 \[
 R_\pm = \{p\in S_+: \star\circ\mathbf{r}(p)=1\}.
 \]
For $N\in\ZZ_+$, let $R_{\pm,N}= \{p\in R_{\pm}: p<N\}$ and  let $S_{\pm,N}= \{p\in S_{\pm}: p<N\}$. 

Recall the map $\star: \MKG \to \pm 1 $ from Definition \ref{StarMG}\switch{StarMG} and recall the map $\mathbf{N}: \MKG \to \pm 1$ from Definition \ref{defn:NMG}\switch{defn:NMG}. Both were first shown to be well defined out of $\MK$ and in this proof we use the notation $\star$ and $\mathbf{N}$ to refer to the corresponding maps from $\MK \to \pm 1$, noting that both $\star$ and $\mathbf{N}$ commute with the Galois action so they are each constant on any fixed Galois orbit.
Let $\star_+$ denote the restriction of $\star$ to 
\[\MK^+\defeq \{\alpha \in \MK: \mathbf{N}(\alpha)=1\}\]
and let $\star_-$ denote the restriction of $\star$ to 
\[\MK^-\defeq \{\alpha \in \MK: \mathbf{N}(\alpha)=-1\}.\]

We will prove that 
\begin{equation}\label{eqn:pmdensityker}
d(R_+|S_+) = \frac{\# \ker(\star_+)}{ \# \MK^+}
\quad \text{and} \quad
d(R_-|S_-) = \frac{\# \ker(\star_-)}{ \# \MK^-}.
\end{equation}

Then since $K$ is cyclic of odd degree and 2 is inert in $K/\QQ$, we can apply Lemma 3.3\switch{(Mundy)} in \cite{APSR} to get that $\#\MK = 2^n$. Note that $\MK= \MK^+ \sqcup \MK^-$, a disjoint union. Then $\#\MK^+=\#\MK^- = 2^{n-1}$. Since $n$ is prime all Galois orbits of $\MK$ have size 1 or $n$. By Lemma 3.3\switch{(Mundy)} in \cite{APSR}, the only orbits of size 1 are $\pm 1$. By Lemma 5.2\switch{(pm1)} in \cite{APSR}, $\star(1)=1$ and $\star(-1)=-1$ so $\#\ker(\star_+) = 1 +m_+n$ and $\#\ker(\star_-) = m_-n$.

We now show equation (\ref{eqn:pmdensityker}).
Let $\pm$ denote a fixed sign and let $\mp$ denote the opposite sign. Let $S_{\pm,N}'$ denote the set of primes of $K$ laying above primes in $S_{\pm,N}$ and let $R_{\pm,N}'$ denote the set of primes of $K$ laying above primes in $R_{\pm,N}$.
Since primes in $S$ split completely, 
\[
\frac{\#R_{\pm,N}}{\#S_{\pm,N}}= \frac{\#R_{\pm,N}'}{\#S_{\pm,N}'}.
\]
Let $\mathbf{r}_{0,N}$ denote the restriction of $\mathbf{r}_{0}$ to $S_{\pm, N}'$. Then $R_{\pm,N}'$ is the disjoint union
\[
R_{\pm,N}' = \bigsqcup_{\star(\alpha)=1} \mathbf{r}_{0,N}\inv(\alpha)
\]
taken over elements $\alpha\in\MK^\pm$ such that $\star(\alpha)=1$, or equivalently taken over $\alpha\in\MK$ such that $\star(\alpha)=1$ and $\mathbf{N}(\alpha)=\pm 1 \bmod 4$.
Therefore
\[
\frac{\#R_{\pm,N}'}{\#S_{\pm,N}'} = \sum_{\substack{\mathbf{N}(\alpha)=\pm \\ \star(\alpha)=1}} \frac{\# \mathbf{r}_{0,N}\inv(\alpha)}{\#S_{\pm,N}'}
\]
over elements $\alpha\in\MK$ such that $\star(\alpha)=1$ and $\mathbf{N}(\alpha)=\pm 1 \bmod 4$. Therefore if the limits exist,
\begin{align*}
d(R_\pm|S_\pm) = d(R_\pm'|S_\pm') 
&= \sum_{\substack{\mathbf{N}(\alpha)=\pm \\ \star(\alpha)=1}} \lim_{N\to\infty} \frac{\# \mathbf{r}_{0,N}\inv(\alpha)}{\#S_{\pm,N}'}\\
& =  \sum_{\substack{\mathbf{N}(\alpha)=\pm \\ \star(\alpha)=1}} d(\mathbf{r}_{0}\inv(\alpha)\cap S_{\pm}' | S_{\pm}').
\end{align*}
By Lemma \ref{equidistribution:extpm}\switch{equidistribution:extpm}, for all $\alpha\in \MK$, 
\[
d(\mathbf{r}_{0}\inv(\alpha)\cap S_{\pm}' | S_{\pm}') = \frac{1}{\#\MK^\pm} = \frac{1}{2^{n-1}}.
\]
Therefore
\[
d(R_\pm|S_\pm) = \frac{\#\ker(\star_\pm) }{\#\MK^\pm}
\]
which proves equation (\ref{eqn:pmdensityker}), completing the proof.

%the equidistribution lemma, Lemma 4.3 in \cite{APSR} uses Dirichlet's Density Theorem. The equidistribution lemma uses that $K$ is totally real, cyclic, $U_T=U^2$, and $h(K)$ is odd.

\end{proof}
}

%%%%%%%%%%%%%%

%\section{Pretty Picture}%[pmstarlight3]%[prettypicture]%[treasuremap]%[golden2inert]

%%%%%%%%%%%%

\section{The Existence of Non-Trivial Starlight}\label{sec:prettypicture:yellow}% [golden2inert][prettypicture][pmstarlight3]

%\begin{center}
%    \switch{sec:prettypicture:yellow}
%\end{center}

%\subsection*{Reducing Assumptions in the Cubic Case}% [golden2inert]

In this section we prove that the negative starlight invariant, $m_-$ is equal to 1 and the positive starlight invariant, $m_+$ is equal to 0 for cubic cyclic number fields with odd class number containing a golden unit. To do so, we prove that the $\Gal(K/\QQ)$ orbit of $\MK$ corresponding to a golden unit is of non-trivial size and satisfies $\star$.

\begin{prop}\cite{AF}\label{AFUT2}\switch{AFUT2}
Let $K$ be a cyclic cubic number field with odd class number. Then every totally positive unit is a square. 
%Assume 2 and 5 are inert in $K/\QQ$ and the conductor
\end{prop}
\switchpf{
\begin{proof}
%Let $h(K)$ denote the class number of $K$. 
Let $U\defeq  \OK^\times$ denote the group of units, $U_T$ the totally positive units, and $U^2$ the square units. 
Observe $U^2\subseteq U_T \subseteq U$. Then we have a surjective homomorphism
\[
\phi: \frac{U}{U^2} \to \frac{U}{U_T}.
\]
If none of the nontrivial class representatives of $U/U^2$ are totally positive then $\phi$ is injective.  By Theorem $V$ in \cite{AF}, all signatures are represented by units. Square units are always totally positive and there are 8 signatures and 8 classes of units mod squares, so each class of $U/U^2$ must have a different signature. Therefore $U_T$$=$$U^2$.
\end{proof}
}

%We now define a golden unit of a cyclic cubic number field. We then prove in Proposition \ref{prop:prettypicture:existence} that for cyclic cubic number fields with odd class number containing a golden unit, $m_+=0$ and $m_-=1$.

Let $\OK$ denote the ring of integers of $K$, a cyclic cubic number field. Recall the definition of a golden unit.

\goldenunit*

\begin{lem}\label{golden2inert}\switch{golden2inert}
Let $K$ be a cubic cyclic number field containing two units that sum to 1. Then 2 is inert in $K/\QQ$. 
\end{lem}
\switchpf{
\begin{proof}
%%\todo{To Do: See MathOverflow question.} 
%We show that any cubic cyclic number field $K$ with two units summing to $1$ must have the property 2 is inert in $K/\QQ$.

If 2 is not inert in $K/\QQ$ then either 2 is totally ramified %e=3 f=1 g=1
or 2 splits completely. %e=1 f=1 g=3.
In either case, the inertia degree of 2 in the extension $K/\QQ$ is $f=1$ so letting $\pp$ denote a prime of $K$ above $2$, the injection of residue class fields
\[
\operatorname{\iota}: \ZZ/2 \hookrightarrow \OK/\pp
\]
is also a surjection. 

Let $\omega\in\OK^\times$ be a unit. Then $\iota\inv(\omega)\neq 0$. That is, $\iota\inv(\omega)=1$. If in addition $\iota\inv(1-\omega)=1$ then we would have $\iota\inv(1)= \iota\inv(\omega)+\iota\inv(1-\omega) = 0$, a contradiction. Therefore
$\iota\inv(1-\omega)=0$. This implies $2|\Norm(1-\omega)$ so $1-\omega$ is not a unit.
\end{proof}
}

\begin{cor}\label{cor:golden2inert}\switch{cor:golden2inert}
If $K$ is a cyclic cubic number field containing a golden unit, then 2 is inert in $K/\QQ$.
\end{cor}

%%%%%%%%%%%%%%%% prettypicture

%\subsection*{Pretty Picture Proposition}%[prettypicture]%[cor:pmstarlight3]

%\begin{prop}[Pretty Picture]\label{prop:prettypicture:existence}\switch{prop:prettypicture:existence}
%Let $K$ be a cyclic cubic number field with odd class number. Assume $K$ has a golden unit $\omega\in\OK^\times$. %and assume $5$ is inert in $K/\QQ$. 
%Then there exists an orbit $X\in\mathbb{M}_{4,G}$  ($X=[\omega]$) such that
%\begin{description}
 %   \item[i] $\#X>1$,
%    \item[ii] $\star(X)=1$, and
%    \item[iii] $\mathbf{N}(X)=-1$,
%\end{description}
%\begin{align*}
%    &\mathbf{N}(X)=-1, \\
%    &\#X>1, \quad \text{and} \\
%    &\star(X)=1.
%\end{align*}
%\end{prop}

%\begin{defn}\label{defn:goldenunit}
%For $K$ a cyclic cubic number field, define $\omega\in\OK^\times$ to be a golden unit iff
%\[
%\Norm(\omega)=\Norm(1-\omega)=-1.
%\]
%\end{defn}

The next Lemma asserts the existence of a non-trivial orbit of $\MK$ such that $\star$ is true whenever $K$ is a cyclic cubic number field with odd class number containing a golden unit. %such that 5 is inert in $K/\QQ$. 
In particular, for such fields, $\star$ is true for the unique non-trivial orbit in $\MKG$ of negative norm.  %where the action is $\Gal(K/\QQ)$ acting on $\MK$ .

\begin{prop}\label{prop:prettypicture:existence}\switch{prop:prettypicture:existence}
Let $K$ be a cyclic cubic number field with odd class number. 
Assume $K$ has a golden unit $\omega\in\OK^\times$. %and assume $5$ is inert in $K/\QQ$. 
Then there exists an orbit $X\in\mathbb{M}_{4,G}$, $(X=[\omega])$ such that
\begin{description}
    \item[i] $\star(X)=1$, 
    \item[ii] $\mathbf{N}(X)=-1$, and
    \item[iii] $\#X=3$.
\end{description}
%\begin{align*}
%    &\mathbf{N}(X)=-1, \\
%    &\#X>1, \quad \text{and} \\
%    &\star(X)=1.
%\end{align*}
\end{prop}

\begin{comment}
\switch{Proof Ingrediants:
\begin{itemize}
    \item 2 is inert by the existence of a golden unit. See Mathoverflow.
    
    \item $U_T=U^2$ by \cite{AF}.
\end{itemize}
\begin{enumerate}
    \item[i] $\#X>1$ Ingredients:
        \begin{itemize}
            \item Mundy both statements. ($K$ cyclic, odd degree, 2 inert).
            \item Dirichlet's Unit Theorem. ($K$ totally real).
            \item Equidistribution Lemma. ($K$ cyclic, totally real, $U_T=U^2$, odd class number $h(K)$)
            \item 2 inert used directly.
            \item odd class number $h$ used directly.
            \item $\mathbf{x}_K$ well-defined. ( $K$ totally real, Galois over $\QQ$, odd class number $h$, $U_T=U^2$, and 2 is inert in $K/\QQ$. )
        \end{itemize} 
    \item[ii] $\star(X)=1$ Ingredients:
        \begin{itemize}
            \item Definition of a golden unit used directly. ($K$ cyclic, cubic, totally real).
        \end{itemize}
    \item[iii] $\mathbf{N}(X)=-1$ Ingredients:
        \begin{itemize}
            \item Equidistribution Lemma. ($K$ cyclic, totally real, $U_T=U^2$, odd class number $h$).
        \end{itemize}
\end{enumerate}
}
\end{comment}

\switchpf{
\begin{proof}
%%\todo{To Do: re-write this proof more clearly. Make sure everything used is cited or stated earlier.}
%Note that $K$ is totally real and every totally positive unit of $K$ is a square by Lemma \ref{lem:AF:UTU2}.

\begin{description}
\item[i]
Let $\sigma$ be a generator of $\Gal(K/\QQ)$.
Observe that for a golden unit $\omega\in \OK^\times$,
$
\{\omega,\sigma(\omega), \sigma^2(\omega)\} 
$
generates $U/U^2$. 
    We know $\omega \nequiv 1$ in $U/U^2$ because $\Norm(\omega)=-1$. We know $\omega \neq -1$ in $U/U^2$ because otherwise $1-\omega$ would be totally positive but $\Norm(1-\omega)=-1$. Therefore $\omega \nequiv \pm 1$ in $U/U^2$. 
    Suppose (for contradiction) that $\omega\sigma(\omega) \equiv \pm 1$ in $U/U^2$. Since $\Norm(\omega)=-1$, this implies 
    $
    -\sigma^2(\omega) \equiv \pm 1 \text{ in } U/U^2. %\implies \omega \equiv \pm 1 \quad \text{in } U/U^2 
    $
    Applying $\sigma$, this implies $\omega \equiv \pm 1$ in $U/U^2$, a contradiction. Therefore 
\begin{equation}\label{eqn:golden2Mink}
\{\omega,\sigma(\omega), \sigma^2(\omega)\}  \quad \text{generates $U/U^2$ as a $\ZZ/2$-vector space.}
\end{equation}

Since $\omega\in \OK$, $1-\omega\in \OK$ and since $\Norm(1-\omega)=-1$, then $1-\omega$ is a unit. By statement (\ref{eqn:golden2Mink}) there exist some $a_0,a_1,a_2\in\ZZ/2$ such that in $U/U^2$
\[
1-\omega \equiv (\omega)^{a_0} (\sigma(\omega))^{a_1} (\sigma^2(\omega))^{a_2}.
\]

Since $\Norm(1-\omega)=\Norm(\omega)=-1$, we know $a_0+a_1+a_2$ must be odd. Therefore 
\[
\begin{pmatrix}
a_0\\a_1\\a_2
\end{pmatrix} 
\in\left\{
\begin{pmatrix}
1\\1\\1
\end{pmatrix},
\begin{pmatrix}
1\\0\\0
\end{pmatrix},
\begin{pmatrix}
0\\1\\0
\end{pmatrix},
\begin{pmatrix}
0\\0\\1
\end{pmatrix}
\right\}.
\]

Suppose (for contradiction) that 
\[
\begin{pmatrix}
a_0\\a_1\\a_2
\end{pmatrix} 
\in\left\{
\begin{pmatrix}
1\\1\\1
\end{pmatrix},
\begin{pmatrix}
1\\0\\0
\end{pmatrix}
\right\}.
\]

Then in $U/U^2$,
\[
1-\omega \equiv -1 \quad \text{or} \quad 1-\omega \equiv \omega.
\]

%Let 
%\[
%u\defeq  \left\{
%\begin{array}{cc}
%    \omega & \text{if } 1-\omega \equiv -1 \text{ mod squares} \\
%    \frac{1}{\omega} & \text{if } 1-\omega \equiv \omega \text{ mod squares}.
%\end{array}
%\right.
%\]

If $1-\omega \equiv -1$, we can write 
$1-\omega = -v^2$ for some $v\in U$. This implies
$\omega = 1 + v^2$, but $\Norm(\omega)=-1$ and $1+v^2$ is totally positive.
Similarly, if $1-\omega \equiv \omega$, we can write $1-\omega = v^2\omega$ for some $v\in U$, which implies $\frac{1}{\omega}=1+v^2$, a contradiction since $\Norm\left(\frac{1}{\omega}\right)=-1$.

Therefore
\[
\begin{pmatrix}
a_0\\a_1\\a_2
\end{pmatrix} 
\in\left\{
\begin{pmatrix}
0\\1\\0
\end{pmatrix},
\begin{pmatrix}
0\\0\\1
\end{pmatrix}
\right\}.
\]
which implies $1-\omega \equiv \tau(\omega)$ in $U/U^2$ for some generator $\tau$ of $\Gal(K/\QQ)$.

Then applying $\tau^2$ yields
\[
1-\tau^2(\omega) \equiv \omega \quad \text{in } U/U^2.
\]
Therefore there is some unit $z$ such that
\begin{align*}
    &\omega ( 1 - \tau^2(\omega)) = z^2 \\
    \implies & \omega ( 1 +\Norm(w)\tau^2(\omega)) = z^2 \\
    \implies & \omega ( 1 +\omega\tau(\omega)(\tau^2(\omega))^2) = z^2 \\
    \implies & \omega +\tau(\omega)(\omega\tau^2(\omega))^2 = z^2.
\end{align*}
Taking $x\defeq 1$ and $y\defeq \omega\tau^2(\omega)$, we see that
\[
\omega x^2 + \tau(\omega) y^2 = z^2
\]
has a nontrivial solution. That is, $( \omega, \tau(\omega))_2=1$. Therefore  $[\omega]\in\mathbb{M}_{4,G}$ is in the kernel of $\star$.

\item[ii]
 We now show that $\mathbf{N}([\omega])=-1$. By Lemma \ref{equidistribution:extpm}\switch{equidistribution:extpm}, there exists a prime $\pp$ of $K$ such that $\mathbf{r}_0(\pp)=\bar{\omega}$ in $\mathbf{M}_4$. That is, letting $\alpha$ be a totally positive generator of $\pp^{h(K)}$, then $\omega\alpha$ is congruent to a square mod $4\OK$. Let 
 \[
 \omega \alpha = 4\gamma + x^2
 \]
 for some $\gamma,x\in\OK$ with $x$ co-prime to 2. Then since $\Norm(\omega)=-1$,
 \[
 -\Norm(\alpha) = \Norm(4\gamma + x^2).
 \]
 Since $\Norm(4\gamma + x^2) \equiv \Norm(x)^2 \bmod 4\ZZ$, this proves 
 \[
 \Norm(\alpha) \equiv -1 \bmod 4\ZZ.
 \]
 Therefore $\mathbf{N}([\omega])=-1.$

\item[iii]
It remains to show that $\#[\omega]=3$. Consider the natural homomorphism
 \[
 m:U/U^2 \to \MK.
 \]
 
 We will prove that $m$ is an isomorphism. Then since $\omega\nequiv \pm 1$ in $U/U^2$, this implies $m(\omega) \neq m(\pm 1) = \pm 1$. By Lemma 3.3\switch{(Mundy)} in \cite{APSR}, $\pm 1$ are the only invariants of the action of $\Gal(K/\QQ)$ on $\mathbf{M}_4$. Therefore if $m$ is an isomorphism, this shows $\#[\omega]=3$.
 
 Note that Lemma 3.3\switch{(Mundy)} 
 in \cite{APSR}  shows that $\#\MK=2^n$. We also know that $\#U/U^2=2^n$ by Dirichlet's unit Theorem. %Therefore, it suffices to show that $m$ is surjective. 
 Therefore, showing surjectivity of $m$ is sufficient to show injectivity and is therefore sufficient to show that $\#[\omega]>1$.

Let $Y \in \MK$. By %Lemma 4.3 \switch{(equidistribution)} in \cite{APSR}
Lemma \ref{equidistribution:extpm}\switch{equidistribution:extpm}, there exists a prime $\pp$ of $K$ such that $\mathbf{r}_0(\pp)=Y$. That is, there exists a prime $\pp$ of $K$ such that  $\overline{\alpha}=Y$ in $\MK$ where  $\alpha\in\OK$ is a totally positive generator of $\pp^{h(K)}$. %wma \pp splits completely in K/\QQ

Recall the definition of the map $\vK: \PP_K^{2\ell} \to U/U^2$, which was proven to be well-defined in Proposition \ref{treasuremap}\switch{treasuremap}.
Let $u\in U$ be a representative of $\vK(\pp)\in U/U^2$. By definition of $\vK$, this implies $u\alpha$ is a square (element) in the narrow ray class field over $K$ of conductor $\pp\infty$. Therefore 2 is unramified in $L\defeq K(\sqrt{u\alpha})$. Since 2 is unramified in $L/K$, then 
\[
2 \nmid \operatorname{disc}(\mathcal{O}_L/\OK).
\]
Note that \switch{by Chapter 2, pg38 equation (7) of \cite{MilneANT}},
\[
4\beta = \operatorname{disc}(\mathcal{O}_L/\OK) (\mathcal{O}_L:\OK[\sqrt{u\alpha}])^2.
\]
Since $\ord_\pp(4{u\alpha})=h$ is odd, then 
\[
(\mathcal{O}_L:\OK[\sqrt{u\alpha}]) = 2.
\]
Therefore there exist some $a,b\in\OK$ such that
\[
\frac{a+b\sqrt{u\alpha}}{2} \in \mathcal{O}_L.
\]
Therefore,
\begin{align*}
&\Norm_{L/K}\left(\frac{a+b\sqrt{u\alpha}}{2}\right) \in \mathcal{O}_K\\
\implies & a^2 - u\alpha b^2 \equiv 0 \bmod 4\OK\\
\implies & u\alpha \equiv \square \bmod 4\OK.
\end{align*}

Therefore $u\alpha \equiv 1$ in $\MK$ and so $m$ is surjective, completing the proof that $\#[\omega]=3$.

\end{description}
\end{proof}
}

%%%%%%%%%%%%%%[cor:pmstarlight3]
%\newpage
\begin{cor}[$\pm$ Starlight for $n=3$]\label{cor:pmstarlight3}\switch{cor:pmstarlight3}
Let $K$ be a cyclic cubic number field with odd class number. 
Assume $K$ has a golden unit. %$\omega\in\OK^\times$
%and assume $5$ is inert in $K/\QQ$. 
Then $m_- = 1$ and $m_+ = 0$.
\end{cor}

\switchpf{
\begin{proof}
By Proposition \ref{prop:prettypicture:existence}\switch{prop:prettypicture:existence}, $m_- = 1$. Note that by %Mundy's Lemma 
Lemma 3.3\switch{(Mundy)} in \cite{APSR}, there is a unique non-trivial $\Gal(K/\QQ)-$orbit of $\MK$ of each fixed norm, $\pm$. Therefore by Lemma 6.4\switch{(halfstar)} in APSR \cite{APSR}, $m_+=0$.
\end{proof}
}

\section{Main Results}\label{sec:main:Dstarequality3:blue}%[Dstarequality3]

%\begin{center}
  %  \switch{sec:main:Dstarequality3:blue}
%\end{center}

For $S\subseteq R$ sets of rational primes, let 
\[
d(R|S)\defeq  \lim_{N\to\infty} \frac{R_N}{S_N}
\]
where $S_N\defeq \{p\in S: \Norm(p)<N\}$ and similar for $R_N$.

Let $K$ be a {cubic cyclic} number field {with odd class number}.
%%\todo{To Do: Check assumptions on $K$ needed to define $d_\pm$}
%containing a golden unit.
Let $S$ denote the set of odd rational primes that split completely in $K/\QQ$. 
Let 
\begin{align*}
&S_+\defeq \{ p \in S: p \equiv 1 \bmod 4\ZZ \}, \quad \text{and} \\
&S_-\defeq \{ p \in S: p \equiv -1 \bmod 4\ZZ \}.
\end{align*}
For a fixed sign $\pm$, define
\[
 R_\pm \defeq \{p\in S_\pm: \spin(\pp,\sigma) = \spin(\pp,\sigma\inv)  \text{ for all } \sigma\neq 1 \in \Gal(K/\QQ) \}.
 \]
 
\Dstarequality*

\switchpf{
\begin{proof}
The assumptions to Theorem \ref{thm:pmdensityformulas}\switch{thm:pmdensityformulas} are true by Proposition \ref{AFUT2}\switch{AFUT2} and Corollary \ref{cor:golden2inert}\switch{cor:golden2inert}.
Apply Theorem \ref{thm:pmdensityformulas}\switch{thm:pmdensityformulas} and Corollary \ref{cor:pmstarlight3}\switch{cor:pmstarlight3} to get $d_+=d(R_+|S_+) = {1}/{4}$ and $d_-=d(R_-|S_-)= {3}/{4}$.

If the limits exist then
\[
d = d(R|S) = d(R_+|S) + d(R_-|S) \quad \text{and}
\]
\[
d(R_\pm|S) = d(R_\pm|S_\pm)d(S_\pm|S) = \frac{1}{2} d(R_\pm|S_\pm)
\]
by Cebotarev's Theorem (see Theorem VII.13.4 in \cite{NeukirchANT} for Dirichlet density or Theorem 4 in \cite{Serre} for natural density).
Therefore $d=d(R|S)$ is the average of $d_+=d(R_+|S_+)$ and $d_-=d(R_-|S_-)$ given in Theorem \ref{thm:Dstarequality3}.

\end{proof}
}

\section{Acknowledgements}

This research was supported by Cornell University and the Max Planck Institute for Mathematics. The author would also like to thank Ravi Ramakrishna, Brian Hwuang, Christian Maire, Keith Conrad, Carlo Pagano, Sam Mundy, and Gerhard Niklasch for helpful mathematical discussions.

\bibliography{APSRIIbib}

\begin{thebibliography}{FIMR13}

\bibitem[AF67]{AF}
J.~V. Armitage and A.~Fr\"ohlich.
\newblock Classnumbers and unit signatures.
\newblock {\em Mathematika}, 14:94--98, 1967.

\bibitem[FIMR13]{FIMR}
J.~B. Friedlander, H.~Iwaniec, B.~Mazur, and K.~Rubin.
\newblock The spin of prime ideals.
\newblock {\em Invent. Math.}, 193(3):697--749, 2013.

\bibitem[Koc00]{HelmutKoch}
Helmut Koch.
\newblock {\em Number theory}, volume~24 of {\em Graduate Studies in
  Mathematics}.
\newblock American Mathematical Society, Providence, RI, 2000.
\newblock Algebraic numbers and functions, Translated from the 1997 German
  original by David Kramer.

\bibitem[McM18]{APSR}
Christine McMeekin.
\newblock On the asymptotics of a prime spin relation, 2018.
\newblock submitted to the Journal of Number Theory, (in revision),
  arXiv:1807.00892.

\bibitem[Mil08]{MilneANT}
James~S. Milne.
\newblock Algebraic number theory (v3.01), 2008.
\newblock Available at www.jmilne.org/math/.

\bibitem[Mil13]{MilneCFT}
J.S. Milne.
\newblock Class field theory (v4.02), 2013.
\newblock Available at www.jmilne.org/math/.

\bibitem[Neu99]{NeukirchANT}
J\"urgen Neukirch.
\newblock {\em Algebraic number theory}, volume 322 of {\em Grundlehren der
  Mathematischen Wissenschaften [Fundamental Principles of Mathematical
  Sciences]}.
\newblock Springer-Verlag, Berlin, 1999.
\newblock Translated from the 1992 German original and with a note by Norbert
  Schappacher, With a foreword by G. Harder.

\bibitem[Ser81]{Serre}
Jean-Pierre Serre.
\newblock Quelques applications du th\'{e}or\`eme de densit\'{e} de
  {C}hebotarev.
\newblock {\em Inst. Hautes \'{E}tudes Sci. Publ. Math.}, (54):323--401, 1981.

\end{thebibliography}
\bibliographystyle{alpha}

\end{document}